\documentclass{ifacconf}

\usepackage{graphicx}      
\usepackage{natbib}        
\usepackage{amsmath}
\usepackage{subcaption}
\usepackage{xcolor}
\usepackage{mleftright}
\usepackage{mathtools}
\usepackage{tikzit}
\begin{document}
\begin{frontmatter}

\title{Port-Hamiltonian modeling and control of a curling HASEL actuator\thanksref{footnoteinfo}} 

\thanks[footnoteinfo]{This work is supported by the EIPHI Graduate School (contract ANR-17-EURE-0002).}

\author{Nelson Cisneros}, 
\author{Yongxin Wu},
\author{Kanty Rabenorosoa},
\author{Yann Le Gorrec}

\address{FEMTO-ST institute, UBFC, CNRS, Besançon, France.(emails: nelson.cisneros@femto-st.fr, yongxin.wu@femto-st.fr, rkanty@femto-st.fr, yann.le.gorrec@femto-st.fr)}

\begin{abstract}                
This paper is concerned with the modeling and control of a curling Hydraulically Amplified Self-healing Electrostatic (HASEL) actuator using the port-Hamiltonian (PH) approach. For that purpose, we use a modular approach and consider the HASEL actuator as an interconnection of elementary subsystems. Each subsystem is modeled by an electrical component consisting of a capacitor in parallel with an inductor connected through the conservation of volume of the moving liquid to a mechanical structure based on inertia, linear, and torsional springs. The parameters are then identified, and the model is validated on the experimental setup. Position control is achieved by using Interconnection and Damping Assignment-Passivity Based Control (IDA-PBC) with integral action (IA) for disturbance rejection. Simulation results show the efficiency of the proposed controller.
\end{abstract} 


\end{frontmatter}

\section{Introduction}

In recent years, one of the most interesting technologies that have been developed for soft robotic applications is the Hydraulically Amplified Self-healing Electrostatic (HASEL) actuator \citep{doi:10.1126/science.aao6139}.

HASEL actuators blend the advantages of Dielectric Elastomer Actuators (DEAs) and fluid-driven soft actuators, combining the convenience of electrical control, excellent electromechanical performances, extensive design flexibility, and various actuation modes \citep{https://doi.org/10.1002/adma.202003375}. There are different types of  HASEL actuators, such as peano, planar, elastomeric donut,  quadrant donut, high-strain peano, and curling actuators. Some interesting applications of HASEL actuators can be found in the literature: a soft gripper for aerial object manipulation \citep{9494690}, an actuator powering a robotic arm  \citep{doi:10.1126/science.aao6139}, an electro-hydraulic rolling soft wheel \citep{ly2022electro}, a peano actuator for enhanced strain, load, and rotary motion \citep{tian2022peano} and soft-actuated joints based on the hydraulic mechanism used in spider legs \citep{https://doi.org/10.1002/advs.202100916}.

 In this paper, we consider, as a benchmark the design, modeling, and control of a simple curling HASEL actuator that can be used as a basic element in more complex robotic structures such as robotic hands or soft grippers. 
 The considered curling HASEL actuator is based on a strain limiting layer that is added to the traditional HASEL mechanism \citep{doi:10.1126/science.aao6139} to change its motion from linear to angular deformation \citep{https://doi.org/10.1002/adma.202003375}.  It is important to derive a reliable model representing the system's dynamics to control the actuator.  
Many papers dealing with the modeling of HASEL actuators have been recently proposed in the literature.
In \citep{volchko2022model}, Dynamic Mode Decomposition with Control (DMDc) is applied to derive a linear model, approximating the system's dynamics. In \citep{hainsworth2022simulating}, a non-linear reduced-order mass-spring-damper (MSD) model for a linear HASEL actuator is proposed. However, these works do not take into account the non-linear behavior (such as the drift effect) or electrical dynamics of HASEL actuators in the model, which can make the control design more challenging and difficult to implement in a real-world application. 

Port-Hamiltonian (PH) formulations are particularly well adapted to represent multi-physical systems. The PH approach is then an excellent candidate to represent the considered dynamics of the HASEL actuator. Interconnection and Damping Assignment-Passivity Based Control (IDA-PBC) serves as a highly effective tool for generating asymptotically stabilizing controllers for Port-Controlled Hamiltonian (PCH) models. \citep{ortega2002interconnection}. Previous works used PHS to model soft robots with energy shaping and IDA-PBC controllers, showing good results. In \citep{franco2021adaptive} and \citep{franco2021position}, energy shaping controllers are used to control the position of a soft continuum manipulator with a large number of degrees of freedom (DOF). In \citep{9868237}, the IDA-PBC method has been successfully used to control a nonlinear Cosserat rod model using an early lumping approach. More recently, in \citep{9793678}, a PH formulation of a one DOF of HASEL planar actuator with position control using an IDA-PBC with Integral Action (IA) has been proposed. Compare to \citep{9793678} we consider here an actuator with bending motion instead of linear deformation. This introduces nonlinearities
in the interconnection matrix. Furthermore, we aim to capture the end position drift effect.
The main contributions of this paper are:
\begin{itemize}
    \item We modeled a curling HASEL actuator using the port-Hamiltonian approach to capture the actuator's electrical and mechanical dynamics. 
    \item We identified the model compared with experimental data and validated it with different input voltages.
    \item We designed an IDA-PBC controller with integral action to control the end point position of the curling HASEL actuator and reject the input disturbances. 
   
\end{itemize}

This paper is organized as follows: Section 
\ref{ModelcHASEL} presents the experimental setup of the curling HASEL actuator and its modeling under the PH framework. The parameter identification is detailed in Section \ref{modelIdentVal}.  In Section \ref{Control}, the controller design is presented. Section \ref{Simulations} presents the simulation results, and the conclusions are given in Section \ref{Conclusions}.
    
\section{Curling HASEL actuator and its PH modeling}\label{ModelcHASEL}    In this section, we first introduce the experimental curling HASEL actuator that is used as a benchmark. We describe the setup and its working principle with reasonable hypotheses that will be used for modeling purposes. We then derive the PH model for this actuator. 
    \subsection{Experimental setup}\label{expsetup}
\begin{figure}[h]
  \centering
    \includegraphics[width=0.4\textwidth]{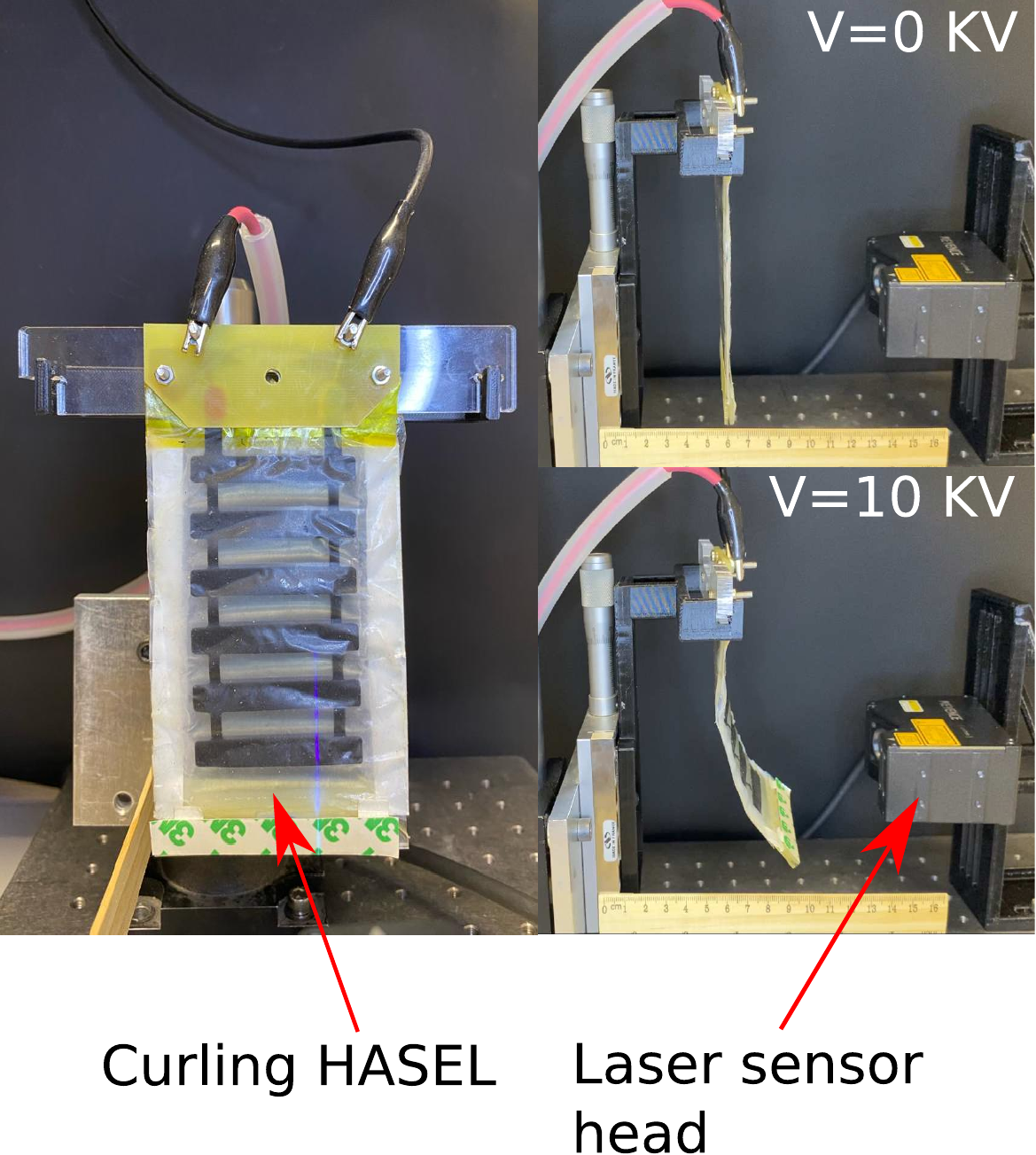}
    \caption{Experimental setup laser sensor and curling HASEL.}
    \label{experimental_setup}
\end{figure}

The experimental setup is shown in Fig. \ref{experimental_setup}. To measure the position, we use a profile laser sensor, Keyence LJ-V7080. We use the high voltage amplifier Trek model 610E. The HASEL actuator used in this work comes from {\it Artimus Robotics}. We use a dSPACE card CLP1104 to receive and send signals from/to the laser position sensor and high voltage amplifier.    
\color{black}
The pictures of Fig. \ref{experimental_setup} (right) show the deformation of the actuator without (upper figure) and with an applied voltage (lower figure). Applying high voltage, the actuator can achieve a horizontal displacement of approximately 3 cm.

\subsection{Curling HASEL actuator description and hypothesis}
 
 The curling HASEL actuator consists of a polymer shell filled with dielectric liquid and half covered by a pair of electrodes attached to a strain-limiting layer to get the bending motion. When an electric field is applied to the electrodes, it creates Maxwell stress acting on the shell that pushes the dielectric liquid inside the shell. This hydraulic pressure changes the shape of the shell and, from the strain-limiting layer, induces the bending of the actuator \citep{https://doi.org/10.1002/adma.202003375}.



\begin{figure}[!htbp]         \centering
         \includegraphics[width=0.18\textwidth, angle= 90]{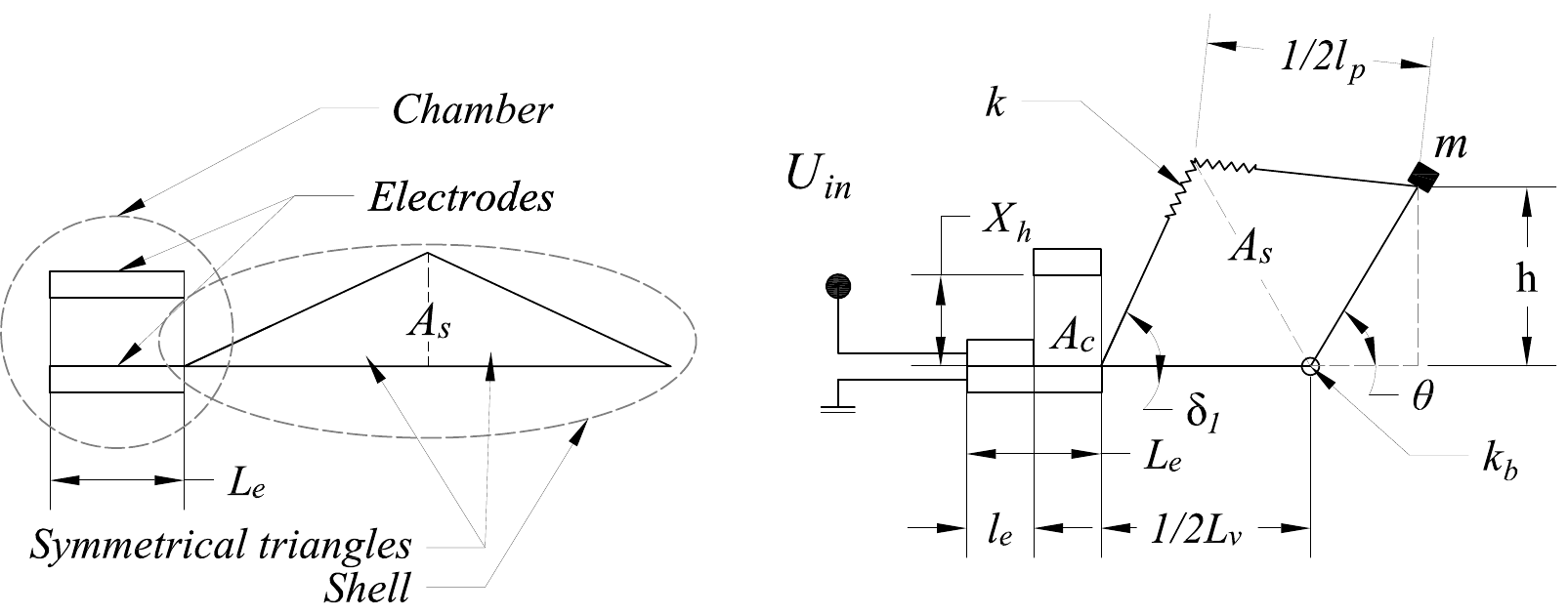}
        \caption{ Basic subsystem. Left: electrodes are totally unzipped. Right: Electrodes are partially zipped when voltage is applied. The shell is deformed.}
        \label{modelschema_1}
\end{figure}

We consider that the actuator depth is uniform, so we proceed to a two-dimensional analysis. We model the actuator using interconnected subsystems. We separate each subsystem model into a chamber and a shell (cf Fig. \ref{modelschema_1}). 

The chamber is the area between the electrodes whilst the shell receives the dielectric liquid when the electrodes are zipped. 
The total volume (the shell's volume plus the chamber's volume) is considered constant, and the dielectric liquid is incompressible. The bending of the bottom film is modeled as a torsional spring.
The top film of the shell is considered to be elongable and contains mechanical energy. The elongation is modeled as a linear spring.
The shell will take on a specific shape based on the volume transferred from the chamber to the shell, resulting from the zipping of the electrodes, which takes place when high voltage is applied.   
The electrodes are modeled as a variable capacitor. The model considers a variable length of the zipped electrodes. The distance between the unzipped electrodes part is considered constant, see Fig. \ref{modelschema_1}.

\subsection{Geometric relations}\label{Geometric_relations}

In this part, we present the relations existing between the angle $\theta$ and the length of the zipped electrodes $l_e$. Then, we derive the actuator position $h(\theta)$ from $\theta$ , i.e., the displacement of the end position of the actuator.

It is crucial to obtain a relation between $\theta$ and $l_e$ because the electrode's capacitance depends on $l_e$. Therefore, it allows us to relate the electrical charge, that depends on the capacitance, with the derivative of the electrical energy respecting the angle $\theta$, joining the electrical and the mechanical part.

We represent the chamber as a rectangular area and the shell is modeled as two symmetric triangles. Fig. \ref{modelschema_1} shows the model variables of a basic subsystem.

The area inside the shell is equal to:
  \begin{equation}\label{eq_A2_0_2}
   A_s=\frac{1}{4}l_p L_v sin (\delta_1)
\end{equation}
with 
\begin{equation}
      \delta_1=\frac{\pi+\theta}{2}-\sin^{-1}\left(\frac{L_v}{l_p}\sin\left(\frac{\pi-\theta}{2}\right)\right).
\end{equation}
The total area is:
\begin{equation}\label{eq_A1_0_2_1}
   {A_{T}=A_s+X_h (L_e-l_e)},
\end{equation}
where $L_v$ and $l_p$ are the lengths of the bottom and the top film. $X_h$ is the height of the chamber and $L_e$ is the length of the chamber.

The zipped electrodes length is then:
\begin{equation}
      l_e=L_e-\frac{1}{X_h}\left(A_T-\frac{L_vl_p}{4}sin (\delta_1)\right). 
\end{equation}

\subsection{Curling HASEL port-Hamiltonian model}\label{model}

This section presents the port-Hamiltonian model for the curling HASEL actuator.
PH formulations consist in using the energy variables as state variables and in writing the dynamics of the system on the form \citep{van2000l2}:
\begin{equation}\label{eqn:phs2}
\begin{array}{ll}
&\dot{x}=[J(x)-R(x)]\frac{\partial{H}}{\partial{x}}(x)+g(x)u;\\
&y=g^T(x)\frac{\partial{H}}{\partial{x}}(x),
\end{array}
\end{equation}
where   $J(x)=-J^T(x)$ is the interconnection matrix, $R(x)=R^T(x)\geq0$ is the dissipation matrix and $H$ is the total energy of the system (Hamiltonian). 

\begin{figure}[ht!]
\centering
\input{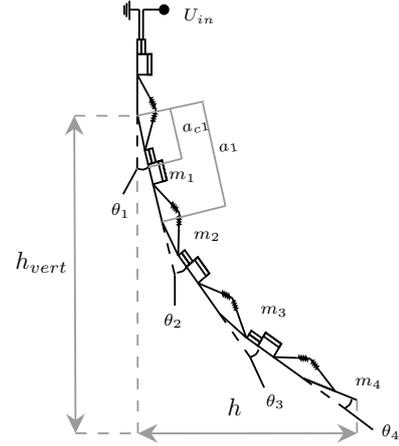}
\caption{Four interconnected subsystems. The same voltage is applied to the entire system.} 
\label{fig:schema4subs_1}
\end{figure}

By combining basic subsystems, we can represent the overall dynamic behavior of the HASEL actuator. In what follows we consider four subsystems (cf Fig. \ref{fig:schema4subs_1}) but the model can be extended to $n \in \mathbb{N}$ subsystems. We consider that the subsystems share the same input voltage. 
 The total energy of the system is given by:
\begin{equation}\label{DE20_2}
\begin{array}{ll}
H(\theta,l_p,p,\phi,Q)=&H_{\theta}(\theta)+H_{l_p}(l_p)+H_{g}(\theta)+\\&H_{p}(p) +H_\phi(\phi,Q){+H_Q(Q,\phi,\theta,l_p)},
\end{array}
\end{equation} 
where
\begin{equation}
   H_{\theta}=\frac{1}{2}
\sum_{i=1}^{n}K_{b_{i}}\theta^2_i=\frac{1}{2}\theta^T K_{b}\theta
\end{equation}
is the potential energy where $K_b=\text{diag}[K_{b_1} \; K_{b_2} \; \ldots \; K_{b_n} ]$  is the stiffness matrix and $\theta=[\theta_1 \; \theta_2  \; \ldots \; \theta_n]$ represents the angular vector of each subsystem,
the second term is the potential energy related to the linear springs:
\begin{equation}
   H_{l_p}=\frac{1}{4}\sum_{i=1}^{n}K_{i}(l_{p_i}-L_{p_i})^2=\frac{1}{4}(l_p-L_p)^T K (l_p-L_p),
\end{equation}
where $K=\text{diag}[K_{1} \; K_{2} \;  \ldots \; K_{n}]$ and $l_p^T=[l_{p_1} \;l_{p_2} \; \ldots \; l_{p_n} ]$,
the third term is the total potential energy related to gravity: 
\begin{align}
   H_{\text{g}}=&\sum_{i=1}^{n}H_{g_i},
\end{align}
and the fourth term the kinetic energy: 
\begin{equation}
H_p=\frac{1}{2}p^TM^{-1} p, 
\end{equation}
where $M$ is the matrix of inertia and $p$ is the vector of angular momentum $p^T=[p_1\; p_2\; \ldots\;p_n]$.

 The electrical energy has two components: the energy related to the magnetic flux that allows us to represent the drift effect and the energy related to the charge. The inductor discharges the capacitor over time.

\begin{align}
    H_{\phi}=\frac{1}{2}\sum_{i=1}^{n}\frac{\phi_i^2}{L_i}=\frac{1}{2}\phi^T L^{-1} \phi 
\end{align}

The energy stored in the capacitor is:
\begin{equation}
\begin{array}{ll}
H_Q&=\frac{1}{2}\sum_{i=1}^{n}\frac{Q_i^2}{C_{s_i}}=\frac{1}{2}Q^TC^{-1}Q,
\end{array}
\end{equation}
where $\phi^T=[\phi_1 \;\phi_2 \; \ldots \; \phi_n]$ is the magnetic flux, $L$ is the inductance of the equivalent electric circuit $L=\text{diag}[L_{1} \;L_{2} \;\ldots \; L_{n}]$. $C=\text{diag} [C_{s_1} \;C_{s_2} \;\ldots \; C_{s_n}]$ is the capacitance of the equivalent electric circuit and $Q=[Q_1 \;Q_2 \; \ldots \; Q_n]^T$ is the charge.
The capacitance of a subsystem is $C_{s_i}=C_{1_i}+C_{2_i}$ where 
$C_{1_i}=\frac{\epsilon_0 \epsilon_r w l_{e_i}}{2t}$ is the capacitance of the zipped part and $C_{2_i}=\frac{\epsilon_0 \epsilon_r w (L_e-l_{e_i})}{2t+X_h}$ the capacitance of the unzipped part. 

The input gain is $ga^T=[ga_1 \; ga_2\; ... \; ga_n]$. To capture the system's nonlinearities, the input gain is a nonlinear function that depends on the angular position $ga_i=\gamma_1 \cos{(\gamma_2 \theta_i)}$. 

The conductance of the equivalent electric circuit is ${\bar{R}}=\text{diag}[\frac{1}{R_1}  \; \frac{1}{R_2} \; ... \; \frac{1}{R_n}]$. The damping coefficient of the system is $b=\text{diag}[b_1\;b_2\; ... \;b_n]$. The resistance associated to the inductance is $r_{L}=\text{diag}[r_{L_1}\;r_{L_2}\; ... \;r_{L_n}]$. We define the term $d=\text{diag}\left(\frac{2A_s}{l_p}\right)$.   
  The proposed port-Hamiltonian model of the curling HASEL actuator is then:
 {\small 
\begin{align}\label{eqModel}
 \underbrace{
 \begin{bmatrix}
\dot{\theta}\\
 \dot{l_{p}}\\ 
  \dot{p}\\ 
  \dot{\phi}\\
  \dot{Q}\\
\end{bmatrix}
}_{\dot{x}}
   =&
\underbrace{   
  \begin{bmatrix}
  0   & 0                        & I                          & 0          & 0           \\
  0   & 0                        & d                          & 0          & 0           \\
 -I   & -d                       & -b                         & 0          & 0           \\
  0   & 0                        & 0                          & -r_{L}     & I           \\
  0   & 0                        & 0                          &-I          & -\bar{R}     \\
\end{bmatrix}
 }_{J-R}
 \underbrace{
  \begin{bmatrix}
  \nabla_{\theta}H\\
  \nabla_{l_p}H\\
  \nabla_{p}H\\
  \nabla_{\phi}H\\
  \nabla_{Q}H
\end{bmatrix}
}_{\nabla_xH}
+
\underbrace{
  \begin{bmatrix}
 0 \\
 0 \\
 0 \\
 0 \\
\bar{R}ga(\theta)
\end{bmatrix}
}_{g}
 U_{in};\\
 y=& \underbrace{(\bar{R}ga(\theta))^TC^{-1}Q}_{g^T\nabla_xH}.\nonumber
\end{align} }

 The output $y=i_e$ is the current that is power-conjugated to the input voltage. 
The energy balance equation can be written as:
 \begin{equation}
 \begin{array}{ll}
     \frac{\partial H}{\partial t} &=-\frac{\partial H^T}{\partial x} R \frac{\partial H}{\partial x} + y^Tu;\\
     \frac{\partial H}{\partial t} &\leq y^Tu =i_e U_{in}.
     \end{array}
 \end{equation} 

The displacement for  $n$ interconnected subsystems is computed as 
\begin{align}\label{h}
    h(\theta)=(L_v+L_e)\left(\sum_{i=1}^{n}\sin(\sum_{j=1}^{i}\theta_j)\right)
\end{align}


\section{Model identification and validation}
\label{modelIdentVal}

\begin{figure}[!htbp]
     \centering
     \begin{subfigure}[b]{0.3\textwidth}
         \centering
         \includegraphics[width=0.8\textwidth, angle= 0]{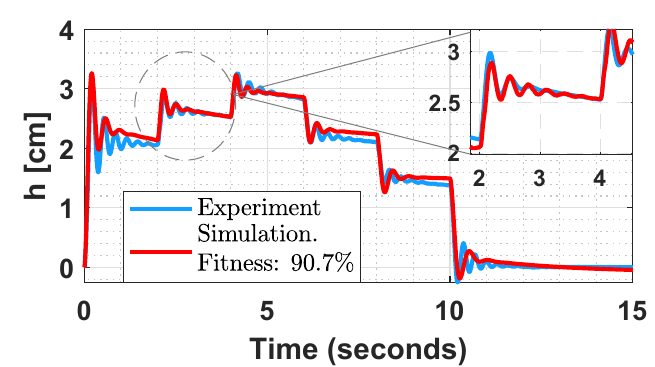}
         \caption{}
         \label{modelresults_1output}
     \end{subfigure}
     \hfill
     \begin{subfigure}[b]{0.17\textwidth}
         \centering
         \includegraphics[width=0.8\textwidth, angle= 0]{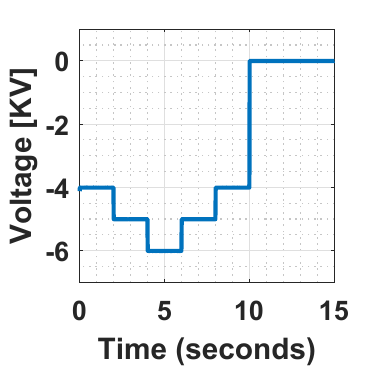}
         \caption{}
         \label{modelresults_1input}
     \end{subfigure}
        \caption{(\ref{modelresults_1output}) Model identification, fitness: 90.7\%. (\ref{modelresults_1input}) Input signal.
      }
        \label{modelresults_1}
\end{figure}

\begin{figure}[!htbp]
     \centering
     \begin{subfigure}[b]{0.31\textwidth}
         \centering
         \includegraphics[width=0.8\textwidth, angle= 0]{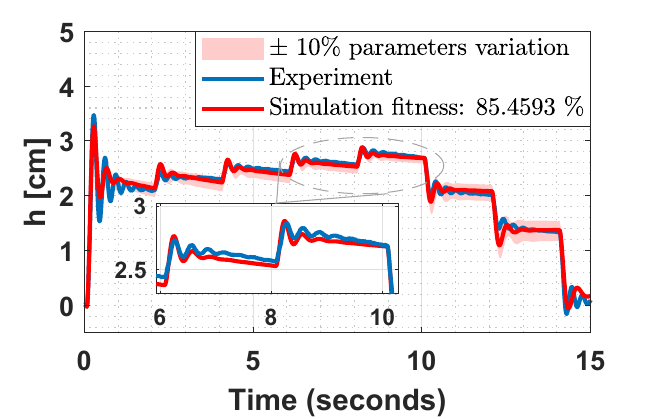}
         \caption{}
         \label{modelresults_2output}
     \end{subfigure}
     \hfill
     \begin{subfigure}[b]{0.17\textwidth}
         \centering
      \includegraphics[width=0.8\textwidth, angle= 0]{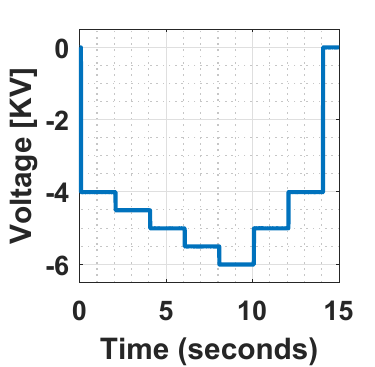}
         \caption{}
         \label{modelresults_2input}
     \end{subfigure}
      \begin{subfigure}[b]{0.31\textwidth}
         \centering
         \includegraphics[width=0.8\textwidth, angle= 0]{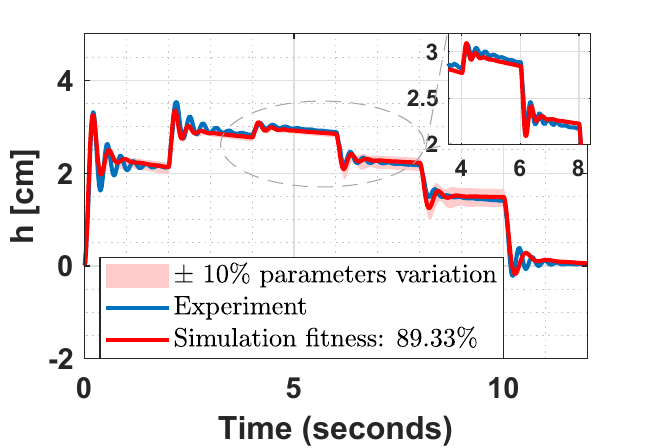}
         \caption{}
         \label{modelresults_4output}
     \end{subfigure}
     \hfill
     \begin{subfigure}[b]{0.17\textwidth}
         \centering
         \includegraphics[width=0.8\textwidth, angle= 0]{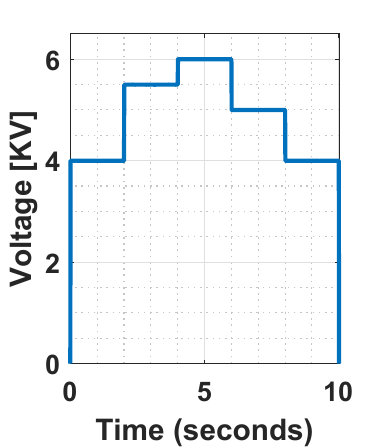}
         \caption{}
         \label{modelresults_4input}
     \end{subfigure}
        \caption{(\ref{modelresults_2output}) Model validation, negative input fitness: 85.46\% (\ref{modelresults_2input}) Input signal. \ref{modelresults_4output}) Model validation, positive input fitness: 89.33\% (\ref{modelresults_4input}) Positive input signal. We can observe the model's behavior in response to a variation of 10\% around the nominal values.
      }
        \label{modelresults_2}
\end{figure}

We identify the key parameters of the system using the experimental data obtained from the experimental setup of Fig. \ref{expsetup}. The Levenberg–Marquardt algorithm is used to find the parameters $K_b$, $b$, $L$, $\gamma_1$ and $\gamma_2$. The identification results are shown in Fig. \ref{modelresults_1} with a fitting of $90.7\%$. Two datasets are used to validate the obtained parameters, one with negative inputs and another with positive inputs, as shown in Fig. \ref{modelresults_2}.

 The fitting between the model and the experimental data is computed using the normalized root mean squared error (NRMSE):

\begin{align}
    fit(i)=\frac{\parallel x_{ref}(:) - x_{data}(:) \parallel}{\parallel x_{ref}(:) - (x_{ref}(:)) \parallel}
\end{align}

 where $\parallel.\parallel$ is the 2-norm of a vector.

One can see the identified model can cope with the main dynamics of the considered Curling HASEL actuator with a fitting comprised between $85\%$ and $89.33\%$. 
 The identified parameters are listed in Table \ref{parameters}.
 \begin{table}[h]
\begin{center}
\resizebox{0.5\textwidth}{!} {
\begin{tabular}{  c  c  c   l  }
\hline
Symbol   &  Value & Units   & Definition  \\
\hline
$L_p$   &0.015   & $m$       & Length of top film \\
$L_v$   &0.015   & $m$       & Length of bottom film \\
$L_e$   &0.015   & $m$       & Length of electrodes \\
$X_h$   & 0.002  &  $m$    & Chamber high       \\
$m$     & 0.047 & $kg$      & Mass     \\ 
$\epsilon_r$         & 2.2 &  $F/m$    & Relative permittivity    \\ 
$\epsilon_0$         & 8.854x$10^{-12}$ &  $F/m$   & Vacuum permittivity    \\ 
$w$      & 0.05&  $m$      & Actuator width    \\ 
$t$      &$18$x$10^-6$  &  $m$      & Film thickness    \\ 
$R_i$      & 10 &  $\Omega$ & Resistance    \\ 
$r_L$      & 20 &  $\Omega$ & Resistance    \\ 
$L$     & 150 & $F$      & Inductance     \\ 
$K$    & 400 &  $ N /m$  & Spring constant    \\ 
$K_b$    & 0.202 &  $ N m/ rad$   & Torsional spring constant     \\ 
 $b$      & 0.0199 & $kgs$     & Damping     \\ 
  $\gamma_1$      & 104.33 & -     & Gain parameter    \\ 
  $\gamma_2$      & 7.67 & -    & Gain parameter    \\
\hline
\end{tabular}
}
\caption{Model parameters.}
\label{parameters}
\end{center}
\end{table}
 

\section{Position control design}\label{Control}

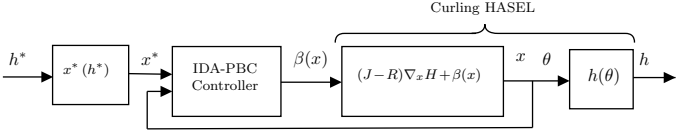
\begin{figure}[ht!]
  \centering
    \resizebox{0.5\textwidth}{!} {

\tikzset{every picture/.style={line width=0.75pt}} 

\begin{tikzpicture}[x=0.75pt,y=0.75pt,yscale=-1,xscale=1]

\draw   (130,54) -- (212,54) -- (212,107.14) -- (130,107.14) -- cycle ;
\draw   (258,54) -- (380,54) -- (380,107.14) -- (258,107.14) -- cycle ;
\draw    (212,81.14) -- (255,81.14) ;
\draw [shift={(258,81.14)}, rotate = 180] [fill={rgb, 255:red, 0; green, 0; blue, 0 }  ][line width=0.08]  [draw opacity=0] (8.93,-4.29) -- (0,0) -- (8.93,4.29) -- cycle    ;
\draw    (380,80.14) -- (427,80.14) ;
\draw [shift={(430,80.14)}, rotate = 180] [fill={rgb, 255:red, 0; green, 0; blue, 0 }  ][line width=0.08]  [draw opacity=0] (8.93,-4.29) -- (0,0) -- (8.93,4.29) -- cycle    ;
\draw    (402,80.14) -- (402,115.64) ;
\draw    (111,116.64) -- (402,115.64) ;
\draw    (111,88.14) -- (127,88.14) ;
\draw [shift={(130,88.14)}, rotate = 180] [fill={rgb, 255:red, 0; green, 0; blue, 0 }  ][line width=0.08]  [draw opacity=0] (8.93,-4.29) -- (0,0) -- (8.93,4.29) -- cycle    ;
\draw    (111,88.14) -- (111,116.64) ;
\draw   (430,56.35) -- (478,56.35) -- (478,101.35) -- (430,101.35) -- cycle ;
\draw    (478,77.6) -- (510,77.6) ;
\draw [shift={(510,77.6)}, rotate = 180] [fill={rgb, 255:red, 0; green, 0; blue, 0 }  ][line width=0.08]  [draw opacity=0] (8.93,-4.29) -- (0,0) -- (8.93,4.29) -- cycle    ;
\draw    (97,74.6) -- (127,75.09) ;
\draw [shift={(130,75.14)}, rotate = 180.94] [fill={rgb, 255:red, 0; green, 0; blue, 0 }  ][line width=0.08]  [draw opacity=0] (8.93,-4.29) -- (0,0) -- (8.93,4.29) -- cycle    ;
\draw   (39.55,53.14) -- (97,53.14) -- (97,98.14) -- (39.55,98.14) -- cycle ;
\draw    (2,77.14) -- (37,77.14) ;
\draw [shift={(40,77.14)}, rotate = 180] [fill={rgb, 255:red, 0; green, 0; blue, 0 }  ][line width=0.08]  [draw opacity=0] (8.93,-4.29) -- (0,0) -- (8.93,4.29) -- cycle    ;
\draw   (479,51) .. controls (479,46.33) and (476.67,44) .. (472,44) -- (376,44) .. controls (369.33,44) and (366,41.67) .. (366,37) .. controls (366,41.67) and (362.67,44) .. (356,44)(359,44) -- (260,44) .. controls (255.33,44) and (253,46.33) .. (253,51) ;

\draw (140,65) node [anchor=north west][inner sep=0.75pt]  [font=\small] [align=left] {\begin{minipage}[lt]{42.51pt}\setlength\topsep{0pt}
\begin{center}
IDA-PBC\\Controller
\end{center}

\end{minipage}};
\draw (220,56) node [anchor=north west][inner sep=0.75pt]   [align=left] {$\displaystyle \beta ( x)$};
\draw (388,58) node [anchor=north west][inner sep=0.75pt]   [align=left] {$\displaystyle x$};
\draw (452,66) node [anchor=north west][inner sep=0.75pt]   [align=left] {$ $};
\draw (439,68.35) node [anchor=north west][inner sep=0.75pt]  [font=\normalsize] [align=left] {\begin{minipage}[lt]{23.27pt}\setlength\topsep{0pt}
\begin{center}
$\displaystyle h( \theta )$
\end{center}

\end{minipage}};
\draw (481,57) node [anchor=north west][inner sep=0.75pt]   [align=left] {$\displaystyle h$};
\draw (105,57) node [anchor=north west][inner sep=0.75pt]   [align=left] {$\displaystyle x^{*}$};
\draw (46,65) node [anchor=north west][inner sep=0.75pt]  [font=\small] [align=left] {$\displaystyle x^{*}\left( h^{*}\right)$};
\draw (5,57) node [anchor=north west][inner sep=0.75pt]   [align=left] {$\displaystyle h^{*}$};
\draw (268,69) node [anchor=north west][inner sep=0.75pt]  [font=\footnotesize] [align=left] {\begin{minipage}[lt]{70.37pt}\setlength\topsep{0pt}
\begin{center}
$\displaystyle ( J-R) \nabla _{x} H+\beta ( x)$
\end{center}

\end{minipage}};
\draw (408,57) node [anchor=north west][inner sep=0.75pt]   [align=left] {$\displaystyle \theta $};
\draw (320,19) node [anchor=north west][inner sep=0.75pt]  [font=\small] [align=left] {\begin{minipage}[lt]{64.45pt}\setlength\topsep{0pt}
\begin{center}
Curling HASEL
\end{center}

\end{minipage}};

\end{tikzpicture}

}
  \caption{Closed-loop scheme with the controller $\beta(x)$ inputs are the state variables $x$ and the desired values $x^*$. The system input is the necessary voltage computed by the controller. $h(x)$ (\ref{h}) is the function that allows us to find the final position as a function of each link angle.} 
  \label{controlschema_1}
\end{figure}

In this work, we aim to control the endpoint position of the curling HASEL actuator (denoted by $h$). To this end, we propose an IDA-PBC design method. This method aims to find a state feedback control law $\beta(x)$ to map the open-loop system to a desired closed-loop system of the form:  
\begin{equation}
    \dot{x}=(J_d-R_d)\nabla_xH_d 
\end{equation}
with the desired interconnection and damping matrices $J_d$, $R_d$ and the desired energy function $H_d$ in the closed-loop system.
The control scheme is shown in Fig. \ref{controlschema_1}. 

The desired equilibrium points $x^* = [\theta^*,l_p^*,0,\phi^*, Q^*]^T$ which can be computed from the desired endpoint position $h^*$ (solving (\ref{h})) and the state variables  $x = [\theta,l_p,p,\phi,Q]^T$ are the controller ($\beta(x)$) inputs. 
We define as desired interconnection and dissipation matrices:
    \begin{equation}\label{Jd}
   J_d - R_d
=
\begin{bmatrix}
  0         &  0         &  J_{13}            &             0  &\alpha_{1}   \\
  0         &  0         &  J_{23}            &             0  &\alpha_{2}   \\
  -J_{13}   &  -J_{23}   & -r_{33}               &        J_{43}  &\alpha_{3}   \\
  0         &  0         & -J_{43}            &             0  &\alpha_{4}   \\
  -\alpha_{1} &  -\alpha_{2} & -\alpha_{3}    &  -\alpha_{4}   &-r_{55}
\end{bmatrix},
 \end{equation} 
 where $J_{13}$, $J_{23}$, $J_{43}$, $\alpha_1$, $\alpha_2$, $\alpha_3$ and  $\alpha_4$ are the control design parameters to be determined.
The desired closed loop energy function is defined from the desired equilibrium position of the actuator as:
\begin{equation}
\begin{array}{ll}
    H_d=&(\theta-\theta^*)^T \Tilde{K_b}(\theta-\theta^*)+ (l_{p}-l_{p}^*)^T\Tilde{K}(l_{p}-l_{p}^*)\\&+ p^T M^{-1} p  +(\phi-\phi^*)^T\Tilde{K_\phi}(\phi-\phi^*)\\&+(Q-Q^*)^T\Tilde{K_Q}(Q-Q^*)
    \end{array}
\end{equation}
 
The derivative of $H_d$ with respect to $x$ is given by:
     \begin{equation}
  \nabla_xH_d
=
\begin{bmatrix}
  \Tilde{K_b}(\theta-\theta^*)      \\
  \Tilde{K}(l_{p}-l_{p}^*)          \\
  M^{-1} p                          \\
  \Tilde{K_{\phi}}(\phi-\phi^*)                \\  
  \Tilde{K_Q}(Q-Q^*)                \\  
\end{bmatrix}.
 \end{equation}
   To get the state feedback matching the closed-loop system with a desired PH system $\dot{x}=(J_d-R_d)\nabla_x H_d$ defined above we need to solve the following matching equation:
  \begin{equation}
  g^\perp [J-R]\nabla_x H=g^\perp [J_d-R_d]\nabla_x H_d,
  \end{equation}
 with $g^\perp$ is a full rank annihilator of the input matrix $g$.
We choose the annihilator as follows:
   \begin{equation}
   g^\perp 
=
\begin{bmatrix}
  1  &  0 & 0 & 0 & 0\\
  0  &  1 & 0 & 0 & 0\\
  0  &  0 & 1 & 0 & 0\\
  0  &  0 & 0 & 1 & 0
\end{bmatrix}.
 \end{equation}

  We find $J_{13}$, $J_{23}$, $J_{43}$ as a function of $\alpha_1$, $\alpha_2$, $\alpha_3$ and $\alpha_4$.
  \small
 \begin{align}
  J_{13} = & \text{diag}((\text{diag}(M^{-1} p))^{-1}(M^{-1} p-\alpha_1 \Tilde{K_Q}(Q-Q^*)));\\
  J_{23} = & \text{diag}((\text{diag}(M^{-1} p))^{-1}(dM^{-1} p-\alpha_2 \Tilde{K_Q}(Q-Q^*)));\\
   r_{33}=& \text{diag}((\text{diag}(M^{-1} p))^{-1}(\nabla_{\theta}H + d \nabla_{l_p}H+b M^{-1} p \\&\nonumber
   +\alpha_3 \Tilde{K_Q}(Q-Q^*) - J_{13} \Tilde{K_b}(\theta-\theta^*)-J_{23} \Tilde{K}(l_p-l^*_p)\\&\nonumber
   +J_{43}\tilde{K_{\phi}}(\phi-\phi^*))).\\
  J_{43} =&  \text{diag}((\text{diag}(M^{-1} p))^{-1}(r_L L^{-1}\phi-C^{-1}Q+\alpha_4 \Tilde{K_Q}(Q-Q^*));
\end{align}
   
\normalsize

We obtain the control law considering the next design parameters $\alpha_1=I$, $\alpha_2=0$, $\alpha_3=I$ and $\alpha_4=0$.

\begin{equation}\label{eq:IDAPBC_control}
\begin{array}{ll}
     \beta(x)= &(\bar{R}ga^T\bar{R}ga)^{-1}\bar{R}ga^T( - \Tilde{K_b}(\theta-\theta^*)  - M^{-1}p\\&
     - r_{55} \Tilde{K_Q}(Q-Q^*)+ L^{-1}\phi+(\bar{R}C^{-1}Q)),
       \end{array}
     \end{equation} 

Given fixed values for $l^*_{p}$ and $\theta^*$ from the desired endpoint position $h^*$, we can determine $Q^*$ from the model at steady state. 
\color{black}

\subsection{Disturbance rejection using Integral Action}\label{DistReject}
\vspace{-.2cm}
In this subsection, we propose to improve the robustness of the controller \eqref{eq:IDAPBC_control} to two types of disturbances acting on the actuator using a structure-preserving integral action. The first one is the unknown mass load, which can be regarded as the unactuated external force disturbance ($d_u$). The other one is the disturbance on the actuated input ($d_a$) i.e. the input voltage. Thus, the disturbed closed-loop system with the previous proposed IDA-PBC control law \eqref{eq:IDAPBC_control} can be written as:
\begin{equation}
\begin{array}{ll}
    \begin{bmatrix}\dot{Q}\\
\dot{\theta}\\
 \dot{l_{p}}\\ 
  \dot{p}\\
  \dot{\phi}
\end{bmatrix}=\begin{bmatrix}J_d-R_d\end{bmatrix}\begin{bmatrix}
 \nabla_Q H_d\\
 \nabla_\theta H_d \\
 \nabla_{l_p} H_d \\
 \nabla_p H_d\\
 \nabla_\phi H_d
\end{bmatrix} + \begin{bmatrix}d_a\\ 0 \\ 0 \\ d_u \\0\end{bmatrix},
\end{array}
\end{equation}
where the desired interconnection and the damping matrix are defined as:
\scriptsize
    
\begin{equation}
  J_d(x):=
  \begin{bmatrix}
    \begin{array}{c|c}
  J_{aa}(x) & J_{au}(x) \\
  \hline
 -J^T_{au}(x) & J_{uu}(x)
    \end{array}
  \end{bmatrix}
  =  \begin{bmatrix}
    \begin{array}{c|c c c c}
  0        & -\alpha_1 & -\alpha_2 & -\alpha_3 & -\alpha_4\\
  \hline
  \alpha_1 & 0         &   0       & -J_{13}   & 0\\
  \alpha_2 & 0         &   0       & -J_{23}   & 0\\
  \alpha_3 & J_{13}    &   J_{23}  & 0         & -J_{43}\\
  \alpha_4 & J_{13}    &   0       & J_{43}    & 0\\
    \end{array}
  \end{bmatrix};
\end{equation}
\begin{equation}
  R_d(x):=
  \begin{bmatrix}
    \begin{array}{c|c}
  R_{aa}(x) & R_{au}(x) \\
  \hline
 R^T_{au}(x) & R_{uu}(x)
    \end{array}
  \end{bmatrix}
  =  \begin{bmatrix}
    \begin{array}{c|c c c c}
  r_{55}          & 0           &   0  & 0    & 0 \\
  \hline
  0                      & 0           &   0  & 0    & 0 \\
  0                      & 0           &   0  & 0    & 0 \\
  0                      & 0           &   0  & r_{33}  & 0 \\
  0                      & 0           &   0  & 0    & 0 \\
    \end{array}
  \end{bmatrix}.
\end{equation}
\normalsize
Using the method described in \citep{ferguson2017new} we choose the new closed-loop Hamiltonian as:
\begin{equation}
    H_{cl} = H_d + \frac{K_{int}}{2}\left(
    Q - x_c\right)^2
\end{equation}
and the new closed loop system can be derived as: 
\begin{equation}
    \begin{bmatrix}
      \dot{x}_a\\\hline 
      \dot{x}_u\\\hline 
      \dot{x}_c
    \end{bmatrix} = 
    (J_{cl}-R_{cl}) \begin{bmatrix}
     \nabla_{ x_a} H_{cl}\\\hline 
     \nabla_{ x_u} H_{cl}\\\hline 
     \nabla_{ x_c} H_{cl}
    \end{bmatrix} + \begin{bmatrix}
    d_a \\\hline  0\\0\\d_u\\0
\\\hline 0
    \end{bmatrix}.
\end{equation}
The structure-preserving Integral Action (IA) controller is then given by:
\begin{equation}
\begin{array}{ll}
     \dot{x_c}=&-R_{c_2}(x) \nabla_{x_a}H+(J_{au}+R_{au})\nabla_{x_u}H,\\
     u_{int}=&[-J_{aa}+R_{aa} + J_{c_1}(x)-R_{c_1}(x)-\\&R_{c_2}(x)]\nabla_{x_a}H+ [J_{c_1}(x)-\\&R_{c_1}(x)]K_{int}(x_a-x_c)+2R_{au}\nabla_{x_u}H;
\end{array}
\end{equation}
where $u_{int}$ is the output of the IA controller. $x_c$ is the IA controller state. The actuated state is the charge $Q_m$ whilst the unactuated states are the angle $\theta$, the length $l_{p}$, the angular momentum $p$, and the magnetic flux $\phi$.

From  $J_{au}$      =  $-[\alpha_1 \ 0 \ \alpha_3 \ 0]$,  $J_{c_1}$  = $R_{c_2}$ = $R_{au}$  = $0$  and 
$R_{c_1}$=  $r_{55}$ we obtain the control law: 
\begin{equation}\label{eq:IAcontrol}
\begin{array}{ll}
          u_{int}=&-r_{55}K_{int}(Q-x_c);\\
    \dot{x}_c=&-( \alpha_1  \Tilde{K_b}(\theta-\theta^*)  + \alpha_3  M^{-1} p).
\end{array}
\end{equation}
\vspace{-.2cm}
where the design parameter $K_{int}$ is chosen as a vector of dimensions $1\times n$. The dimension of the controller $\beta(x)$ is 1$\times$1. 
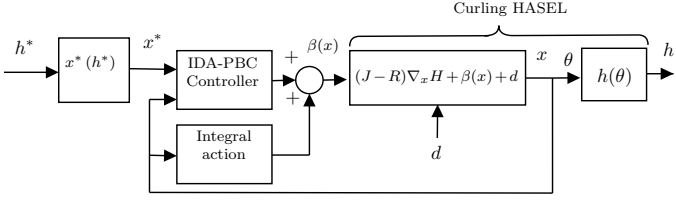
\begin{figure}[!htbp]
  \centering
     \resizebox{0.5\textwidth}{!} {

\tikzset{every picture/.style={line width=0.75pt}} 

\begin{tikzpicture}[x=0.75pt,y=0.75pt,yscale=-1,xscale=1]

\draw   (132,41) -- (198,41) -- (198,82) -- (132,82) -- cycle ;
\draw   (252,41) -- (374,41) -- (374,79) -- (252,79) -- cycle ;
\draw    (99,53.6) -- (129,54.09) ;
\draw [shift={(132,54.14)}, rotate = 180.94] [fill={rgb, 255:red, 0; green, 0; blue, 0 }  ][line width=0.08]  [draw opacity=0] (8.93,-4.29) -- (0,0) -- (8.93,4.29) -- cycle    ;
\draw    (374,60.14) -- (392.72,60.14) -- (409,60.02) ;
\draw [shift={(412,60)}, rotate = 179.58] [fill={rgb, 255:red, 0; green, 0; blue, 0 }  ][line width=0.08]  [draw opacity=0] (8.93,-4.29) -- (0,0) -- (8.93,4.29) -- cycle    ;
\draw    (392.72,60.14) -- (392.73,140.6) ;
\draw    (113,140.14) -- (392.73,140.6) ;
\draw    (113,75.14) -- (129,75.14) ;
\draw [shift={(132,75.14)}, rotate = 180] [fill={rgb, 255:red, 0; green, 0; blue, 0 }  ][line width=0.08]  [draw opacity=0] (8.93,-4.29) -- (0,0) -- (8.93,4.29) -- cycle    ;
\draw    (113,75.14) -- (113,140.14) ;
\draw   (132,93) -- (198,93) -- (198,134) -- (132,134) -- cycle ;
\draw    (113,114.14) -- (129,114.53) ;
\draw [shift={(132,114.6)}, rotate = 181.38] [fill={rgb, 255:red, 0; green, 0; blue, 0 }  ][line width=0.08]  [draw opacity=0] (8.93,-4.29) -- (0,0) -- (8.93,4.29) -- cycle    ;
\draw    (198,114) -- (224,114) ;
\draw    (224,74.1) -- (224,114) ;
\draw [shift={(224,71.1)}, rotate = 90] [fill={rgb, 255:red, 0; green, 0; blue, 0 }  ][line width=0.08]  [draw opacity=0] (8.93,-4.29) -- (0,0) -- (8.93,4.29) -- cycle    ;
\draw   (413,39.35) -- (458,39.35) -- (458,76) -- (413,76) -- cycle ;
\draw    (458,57.14) -- (476.72,57.14) ;
\draw [shift={(476.72,57.14)}, rotate = 180] [fill={rgb, 255:red, 0; green, 0; blue, 0 }  ][line width=0.08]  [draw opacity=0] (8.93,-4.29) -- (0,0) -- (8.93,4.29) -- cycle    ;
\draw    (312.97,81.6) -- (312.76,102.08) ;
\draw [shift={(313,78.6)}, rotate = 90.59] [fill={rgb, 255:red, 0; green, 0; blue, 0 }  ][line width=0.08]  [draw opacity=0] (8.93,-4.29) -- (0,0) -- (8.93,4.29) -- cycle    ;
\draw   (50,32.14) -- (99,32.14) -- (99,77.14) -- (50,77.14) -- cycle ;
\draw    (12,56.14) -- (47,56.14) ;
\draw [shift={(50,56.14)}, rotate = 180] [fill={rgb, 255:red, 0; green, 0; blue, 0 }  ][line width=0.08]  [draw opacity=0] (8.93,-4.29) -- (0,0) -- (8.93,4.29) -- cycle    ;
\draw   (214.5,61.6) .. controls (214.5,56.35) and (218.75,52.1) .. (224,52.1) .. controls (229.25,52.1) and (233.5,56.35) .. (233.5,61.6) .. controls (233.5,66.85) and (229.25,71.1) .. (224,71.1) .. controls (218.75,71.1) and (214.5,66.85) .. (214.5,61.6) -- cycle ;
\draw    (198,62) -- (211.5,61.67) ;
\draw [shift={(214.5,61.6)}, rotate = 178.61] [fill={rgb, 255:red, 0; green, 0; blue, 0 }  ][line width=0.08]  [draw opacity=0] (8.93,-4.29) -- (0,0) -- (8.93,4.29) -- cycle    ;
\draw    (233.5,61.6) -- (247,61.27) ;
\draw [shift={(250,61.2)}, rotate = 178.61] [fill={rgb, 255:red, 0; green, 0; blue, 0 }  ][line width=0.08]  [draw opacity=0] (8.93,-4.29) -- (0,0) -- (8.93,4.29) -- cycle    ;
\draw   (459,37) .. controls (458.98,32.33) and (456.64,30.01) .. (451.97,30.03) -- (365.47,30.46) .. controls (358.8,30.49) and (355.46,28.18) .. (355.43,23.51) .. controls (355.46,28.18) and (352.14,30.52) .. (345.47,30.55)(348.47,30.54) -- (258.97,30.98) .. controls (254.3,31) and (251.98,33.34) .. (252,38.01) ;

\draw (134,44) node [anchor=north west][inner sep=0.75pt]  [font=\small] [align=left] {\begin{minipage}[lt]{42.51pt}\setlength\topsep{0pt}
\begin{center}
IDA-PBC\\Controller
\end{center}

\end{minipage}};
\draw (254,53) node [anchor=north west][inner sep=0.75pt]  [font=\footnotesize] [align=left] {\begin{minipage}[lt]{84.87pt}\setlength\topsep{0pt}
\begin{center}
$\displaystyle ( J-R) \nabla _{x} H+\beta ( x) +d$
\end{center}

\end{minipage}};
\draw (221,31) node [anchor=north west][inner sep=0.75pt]  [font=\small] [align=left] {$\displaystyle \beta ( x)$};
\draw (107,28) node [anchor=north west][inner sep=0.75pt]   [align=left] {$\displaystyle x^{*}$};
\draw (380,40) node [anchor=north west][inner sep=0.75pt]   [align=left] {$\displaystyle x$};
\draw (141,95) node [anchor=north west][inner sep=0.75pt]  [font=\small] [align=left] {\begin{minipage}[lt]{33.34pt}\setlength\topsep{0pt}
\begin{center}
Integral\\action
\end{center}

\end{minipage}};
\draw (446,57) node [anchor=north west][inner sep=0.75pt]   [align=left] {$ $};
\draw (419,51.35) node [anchor=north west][inner sep=0.75pt]  [font=\normalsize] [align=left] {\begin{minipage}[lt]{23.27pt}\setlength\topsep{0pt}
\begin{center}
$\displaystyle h( \theta )$
\end{center}

\end{minipage}};
\draw (468,34) node [anchor=north west][inner sep=0.75pt]   [align=left] {$\displaystyle h$};
\draw (307,106) node [anchor=north west][inner sep=0.75pt]   [align=left] {$\displaystyle d$};
\draw (399,41) node [anchor=north west][inner sep=0.75pt]   [align=left] {$\displaystyle \theta $};
\draw (53,43) node [anchor=north west][inner sep=0.75pt]  [font=\small] [align=left] {$\displaystyle x^{*}\left( h^{*}\right)$};
\draw (18,32) node [anchor=north west][inner sep=0.75pt]   [align=left] {$\displaystyle h^{*}$};
\draw (205,40) node [anchor=north west][inner sep=0.75pt]   [align=left] {$\displaystyle +$};
\draw (206,68.1) node [anchor=north west][inner sep=0.75pt]   [align=left] {$\displaystyle +$};
\draw (319,6) node [anchor=north west][inner sep=0.75pt]  [font=\small] [align=left] {\begin{minipage}[lt]{64.45pt}\setlength\topsep{0pt}
\begin{center}
Curling HASEL
\end{center}

\end{minipage}};

\end{tikzpicture}

}
  \caption{Closed-loop scheme with IA. } 
  \label{controlschema_2}
\end{figure}
The control scheme is shown in Fig. \ref{controlschema_2} and the interconnection and damping matrix are given by:
\begin{equation}
          J_{cl}
         \coloneqq
         \begin{bmatrix}
           0              &  J_{au}        & 0\\
          -J^T_{au}       &  J_{uu}        & 0\\
          0               & 0              & 0\\
         \end{bmatrix}; \;\;
          R_{cl}
         \coloneqq
         \begin{bmatrix}
           r_{55}    & 0              & r_{55}\\
           0                   & r_{33}         & 0\\
          r_{55}    & 0              & r_{55}\\
         \end{bmatrix}.
\end{equation}


\section{Numerical Simulation}\label{Simulations}

\begin{figure}[h]
  \centering
    \includegraphics[width=0.35\textwidth]{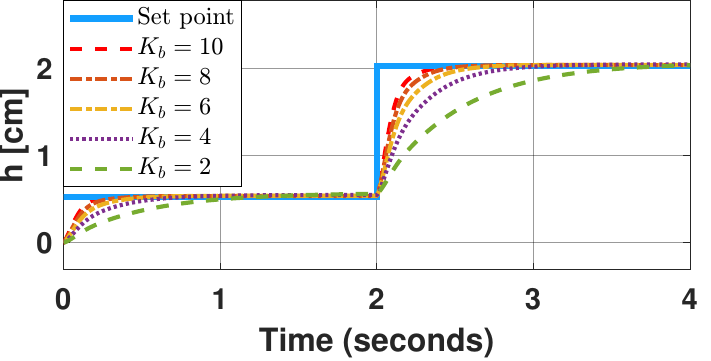}\vspace{-.2cm}
    \caption{Position control  keeping constant the parameters $\Tilde{K_Q}=1000$, the set-point equal to $2$cm and varying the tuning parameter related with the desired angle $\Tilde{K}_b$. }
  \label{control_results_2}
\end{figure}

\begin{figure}
     \centering
     \begin{subfigure}[b]{0.5\textwidth}
         \centering
         \includegraphics[width=0.9\textwidth, angle= 0]{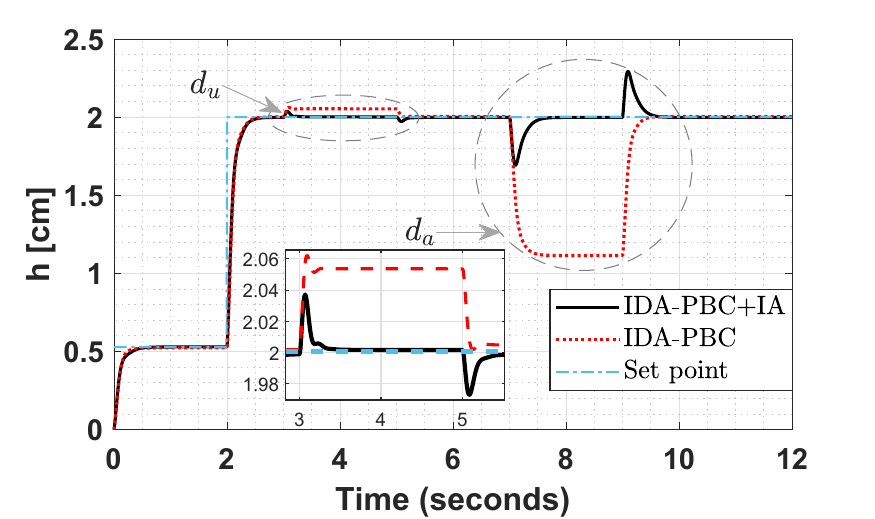}
         \caption{}
         \label{control_results_11}
     \end{subfigure}
     \hfill
     \begin{subfigure}[b]{0.5\textwidth}
         \centering
         \includegraphics[width=0.9\textwidth, angle= 0]{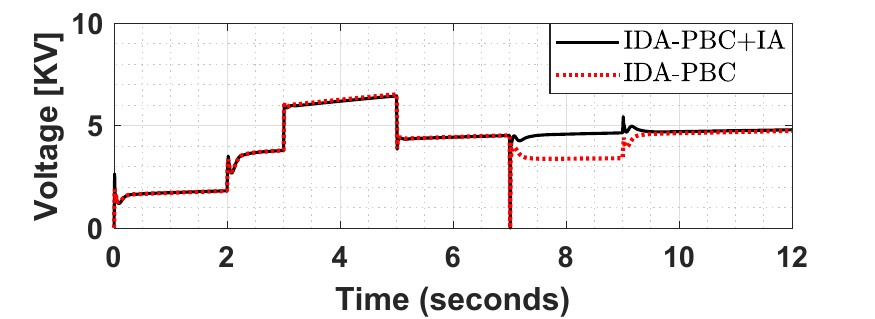}
         \caption{}
         \label{control_results_12}
     \end{subfigure}

         \caption{(\ref{control_results_11}) Position control $\Tilde{K_Q}=1000$, $\Tilde{K_b}=10$ and  $K_{int}=[0.5 \; 11 \; 1.2 \; 0.5]$. The simulation presents disturbances $d_u$ and $d_a$ at $3s$ and $7s$ respectively. (\ref{control_results_12}) IDA-PBC and IDA-PBC+IA control signals.}
    \label{control_results_1}
\end{figure}

In this section, we validate the proposed control methods with numerical simulations. The parameters of the actuator are given in Table. \ref{parameters}.
We implement the IDA-PBC controller \eqref{eq:IDAPBC_control} with the Integral action controller \eqref{eq:IAcontrol} to achieve the desired endpoint position of the curling actuator. To show the different closed loop dynamics performances, we vary the tuning parameter $\Tilde{K_b}$ while keeping the rest of the tuning parameters constant. The controller parameters values are $r_{55}=\text{diag}([0.1 \; 0.1 \; 0.1 \; 0.1])$ and $K_{int}=[0.5 \; 11 \; 1.2 \; 0.5]$. One can observe the endpoint regulation to the desired position in Fig. \ref{control_results_2} using the IDA-PBC method and the rejection of external disturbances in Fig. \ref{control_results_1} with the IDA-PBC+IA method.

Fig. \ref{control_results_2} represents the actuator displacement when $\Tilde{K_b}$ is tuned, while the parameter $\Tilde{K_Q}$ is fixed to a constant value. The desired endpoint position is $h^*(\theta)=2$ cm.  From simulation results shown in Fig. \ref{control_results_2}, the response time of the closed-loop system decreases when $\Tilde{K_b}$ increases, since $\Tilde{K_b}$ can be seen as the actuator's stiffness in the closed-loop system.

Fig. \ref{control_results_1} shows the actuator displacement to different desired set points. The external unactuated disturbance $d_u=-0.04$Nm is added at $3$s during $2$s , and the actuated disturbance $d_a=-30$V is added at $7$s during $2$s (0.3\% of full scale according to the amplifier specifications). The proposed controller with IA can reject the disturbances, while the controller without IA can not reject these disturbances. From the simulation results shown in Fig. \ref{control_results_1}(b), it is seen that the applied controlled voltage on the closed loop actuator always remains below $10$kV, which aligns with the physical consistency of the experimental setup.
\vspace{-0.3cm}

\section{Conclusion}\label{Conclusions}
\vspace{-.3cm}
In this paper we use the port-Hamiltonian framework to model and control a curling HASEL actuator. The actuator's dynamics is divided into two components. The mechanical part of the actuator is characterized by linear and torsional springs, while the deformable capacitor and the inductor represent the electrical part of the system, the coupling being done through the conservation of volume of the overall system. This model can cope with the main dynamic behavior of the actuator with nonlinearities such as the drift effect. An IDA-PBC controller with the integral action is proposed to control the position of the actuator. It is shown that the actuator endpoint position follows the set point and that we can adjust the dynamic performances by varying the tuning parameters of the controller. The use of integral action has improved the robustness of the closed-loop system against external disturbances. 

The perspectives of this work are to implement and validate the proposed IDA-PBC controller with IA in the experimental setup. Furthermore, we intend to model and control more complex structures and soft robots based on HASEL actuators (e.g.,  scorpion, fish, human hands-inspired designs) with the interconnection of basic subsystems.

\small{\bibliography{ifacconf}}             

\end{document}